\documentclass[english]{article}
\usepackage[T1]{fontenc}
\usepackage[latin9]{inputenc}
\usepackage{geometry}
\usepackage{amsmath}
\usepackage{amssymb}
\usepackage{setspace}
\PassOptionsToPackage{normalem}{ulem}
\usepackage{ulem}
\makeatletter
\@ifundefined{date}{}{\date{}}
\makeatother

\usepackage{babel}
\begin{document}
\global\long\def\floor#1{\lfloor#1\rfloor}%

\global\long\def\ceil#1{\lceil#1\rceil}%

\global\long\def\tab{\phantom{tab}}%

\global\long\def\sub{\subseteq}%

\global\long\def\Q{\mathbb{Q}}%

\global\long\def\Z{\mathbb{Z}}%

\global\long\def\id{\text{id}}%

\global\long\def\sgn{\text{sgn}}%

\title{An alternative proof for the irreducibility of the $p$-th cyclotomic
polynomial}
\author{Tom Moshaiov}

\maketitle
Let $p$ be a prime number. As a standard application of the irreducibility
criterion of Eisenstein, it is well known that the $p$-th cyclotomic
polynomial $\Phi_{p}(t)=1+t+\dots+t^{p-1}$ is the minimal polynomial
of $e^{2\pi i/p}$ over $\Q$. This note provides an alternative proof,
utilizing determinants to prove a lemma due to Kronecker.

$ $

$ $

\uline{Claim.} Let $p$ be a prime number, $a_{0},\dots,a_{p-1}$
be indeterminates, and $A=\left(\begin{array}{cccc}
a_{0} & a_{1} & \cdots & a_{p-1}\\
a_{p-1} & a_{0} & a_{1} & \cdots\\
\vdots & \vdots & \ddots & \vdots\\
a_{1} & \cdots & a_{p-1} & a_{0}
\end{array}\right)$ be their circulant matrix. Then $\det A\equiv a_{0}^{p}+a_{1}^{p}+\dots+a_{p-1}^{p}\mod p$.

$ $

$ $

\uline{Lemma (Kronecker).} Let $p$ be a prime number.$\phantom{a}$If$\phantom{.}$
${\displaystyle \sum_{j=0}^{p-1}}a_{j}e^{2\pi ij/p}=0$ with $a_{0},\dots,a_{p-1}\in\Z$,
then $p\mid{\displaystyle \sum_{j=0}^{p-1}}a_{j}$. Equivalently,
if $f(t)\in\Z[t]$ vanishes at some primitive $p$-th root of unity,
then $p\mid f(1)$.

$ $

$ $

\uline{Theorem (Gauss).} Let $p$ be a prime number. Then $\Phi_{p}(t)=1+t+\dots+t^{p-1}$
is the minimal polynomial of $e^{2\pi i/p}$ over $\Q$. Equivalently,
if$\phantom{.}$ ${\displaystyle \sum_{j=0}^{p-1}}a_{j}e^{2\pi ij/p}=0$
with $a_{0},\dots,a_{p-1}\in\Q$ then $a_{0}=a_{1}\dots=a_{p-1}$.

$ $

$ $

\uline{Proof of claim.} Given a permutation $\sigma=(\sigma_{0},\sigma_{1},\dots,\sigma_{p-1})$
of $\Z/p$, we let $T.\sigma=(\sigma_{p-1}+1,\sigma_{0}+1,\sigma_{1}+1,\dots)$,
(additions mod $p$). Geometrically, $T$ translates the matrix rook
arrangement $\{(i,\sigma_{i})\}_{i=0}^{p-1}$ a step down and to the
right, in cyclic fashion. We observe that $T^{p}=\id$, and that the
fixed points of $T$ are the cyclic permutations $(i,i+1,\dots)$,
which produce the products $a_{i}^{p}$. The remaining permutations
are partitioned into $p$-cycles $\sigma,T\sigma,\dots,T^{p-1}\sigma$
all of which produce the same matrix product. We may assume $p$ is
odd, so that $(i,i+1,\dots)$ is even. It remains to show that $\sgn(\sigma)=\sgn(T.\sigma)$.
And indeed, by cyclically rotating $\sigma$, we may assume that $\sigma_{p-1}=p-1$.
Then $\sgn(\sigma)$ and $\sgn(T.\sigma)$ both count the parity of
the number of order inversions in $\sigma_{0},\sigma_{1},\dots,\sigma_{p-2}$.

$ $

$ $

\uline{Proof of Lemma.} Write $\zeta=e^{2\pi i/p}$. The identity
implies $\left(\begin{array}{cccc}
a_{0} & a_{1} & \cdots & a_{p-1}\\
a_{p-1} & a_{0} & a_{1} & \cdots\\
\vdots & \vdots & \ddots & \vdots\\
a_{1} & \cdots & a_{p-1} & a_{0}
\end{array}\right)\left(\begin{array}{c}
1\\
\zeta\\
\vdots\\
\zeta^{p-1}
\end{array}\right)=0$.

The above claim and Fermat's little theorem yield the desired$\phantom{.}$
${\displaystyle \sum_{j=0}^{p-1}}a_{j}\equiv{\displaystyle \sum_{j=0}^{p-1}}a_{j}^{p}\equiv\det(a_{j-i})_{i,j\in\Z/p}=0\mod p$.

$ $

$ $

\uline{Proof of theorem.} The complex roots of $\Phi_{p}(t)$ are
precisely the primitive $p$-th roots of unity. We cannot have

$\Phi_{p}(t)=f(t)g(t)$ for some non-constant $f(t),g(t)\in\Z[t]$,
since $p^{2}\mid f(1)g(1)=\Phi_{p}(1)=p$ would be implied by the
above lemma. Thus $\Phi_{p}(t)$ is irreducible over $\Z[t]$ and
by Gauss\textquoteright s lemma over $\Q[t]$ as well. It follows
that $\Phi_{p}(t)$ is the minimal polynomial of $e^{2\pi i/p}$ over
$\Q$, as we wanted to show.

$ $

$ $

\end{document}